\documentclass{article}
\usepackage{amsmath,amssymb,amsthm, amscd, multicol, mathrsfs, tikz, graphics, float}

\newtheorem{definition}[subsection]{Definition}
\newtheorem{theorem}[subsection]{Theorem}
\newtheorem{corollary}[subsection]{Corollary}
\newtheorem{lemma}[subsection]{Lemma}
\newtheorem{proposition}[subsection]{Proposition}

\newtheorem{observation}[subsection]{Observation}
\renewcommand{\qed}{\hfill$\Box$\\}

\begin{document}
\title{Maximal failed zero forcing sets for products of two graphs}
\author{Prince Allan B. Pelayo\thanks{%
Institute of Mathematics, University of the Philippines, Diliman, Quezon
City 1101; \texttt{papelayo@math.upd.edu.ph}} \; and Ma. Nerissa M. Abara \thanks{%
Institute of Mathematics, University of the Philippines, Diliman, Quezon
City 1101; \texttt{issa@math.upd.edu.ph}}}
\maketitle

\begin{abstract}

Let $G$ be a simple, finite graph with vertex set $V(G)$ and edge set $E(G)$, where each vertex is either colored blue or white. Define the standard zero forcing process on $G$ with the following color-change rule: let $S$ be the set of all initially blue vertices of $G$ and let $u \in S$. If $v$ is the unique white vertex adjacent to $u$ in $G$, color $v$ blue and update $S$ by adding $v$ to $S$. If $S = V(G)$ after a finite number of iterations of the color-change rule, we say that $S$ is a zero forcing set for $G$. Otherwise, we say that $S$ is a failed zero forcing set. In this paper, we construct maximal failed zero forcing sets for graph products such as Cartesian products, strong products, lexicographic products, and coronas. In particular, we consider products of two paths, two cycles, and two complete graphs. \newline
\bigskip


\medskip

 \end{abstract}

\vspace{0.3in}

\section{Introduction}\label{sec1}

A \textit{graph} is an ordered pair $G = (V, E)$ comprised of a non-empty set $V$ of vertices and a (possibly empty) set $E$ of edges. We say that $G$ is a \textit{simple, finite, and undirected graph} if it satisfies all of the following: (a) each edge is a pair of two distinct vertices (that is, $G$ has no loops); (b) each pair of vertices comprises at most one edge (that is, $G$ has no multiple edges); (c) $V$ is finite and; (d) the edges are unordered pairs of vertices. \\
\indent Throughout this paper, we consider simple, finite, and undirected graphs. For terminology and notation not defined here, we refer to \cite{W}. We denote the \textit{vertex set} and \textit{edge set} of a graph $G$ by $V(G)$ and $E(G)$, respectively. The cardinality of $V(G)$ is called the \textit{order of $G$}, denoted $|G|$. An vertex $v$ adjacent to a vertex $u$ is said to be a \textit{neighbor of $u$}. A vertex without any neighbor is called an \textit{isolated vertex}. A set $K \subset V(G)$ of $k$ vertices having the same neighbors outside $K$ is called a \textit{module of order $k$}. Finally, we use the following notations in graph theory:
\begin{center}
\begin{tabular}{lll}
$N(u)$ & \phantom{white} & set of neighbors of vertex $u$ \\
$deg(u)$ & \phantom{white} & degree of vertex $u$ \\
$\Delta(G)$ & \phantom{white} & maximum degree of graph $G$ \\
$\delta(G)$ & \phantom{white} & minimum degree of graph $G$ \\
$G[S]$ & \phantom{white} & subgraph of $G$ induced by set $S$ \\
$c(G)$ & \phantom{white} & number of connected components of $G$\\
$P_n$  & \phantom{white} & path with $n$ vertices \\
$C_n$  & \phantom{white} & cycle with $n$ vertices \\
$K_n$  & \phantom{white} & complete graph with $n$ vertices \\
$G$ $\square$ $H$ & \phantom{white} & Cartesian product of graphs $G$ and $H$ \\
$G$ $\boxtimes$ $H$ & \phantom{white} & strong product of graphs $G$ and $H$ \\
$G$ $\cdot$ $H$ & \phantom{white} & lexicographic product of grapsh $G$ and $H$ \\
$G$ $\circ$ $H$ & \phantom{white} & corona of graph $G$ with $H$ 
\end{tabular}
\end{center}

\indent In this paper, we consider the following graph operations.
\begin{itemize}
\item The \textit{Cartesian product} of two graphs $G$ and $H$, denoted $G$ $\square$ $H$, is the graph with vertex set $V(G)$ $\times$ $V(H)$ such that $(u,v)$ is adjacent to $(u',v')$ if and only if $(1)$ $u=u'$ and $vv'$ $\in$ $E(H)$, or $(2)$ $v=v'$ and $uu'$ $\in$ $E(G)$. In particular, the Cartesian product $P_n$ $\square$ $P_n$ is called a \textit{square grid graph}.
\item The \textit{strong product} of two graphs $G$ and $H$, denoted $G$ $\boxtimes$ $H$, is the graph with vertex set $V(G)$ $\times$ $V(H)$ such that $(u,v)$ is adjacent to $(u',v')$ if and only $(1)$ $uu'$ $\in$ $E(G)$ and $vv'$ $\in$ $E(H)$, or $(2)$ $u=u'$ and $vv'$ $\in$ $E(H)$, or $(3)$ $v=v'$ and $uu'$ $\in$ $E(G)$.
\item The \textit{lexicographic product} of two graphs $G$ and $H$, denoted $G$ $\cdot$ $H$, is the graph with vertex set $V(G)$ $\times$ $V(H)$ such that $(u,v)$ is adjacent to $(u',v')$ if and only $(1)$ $uu'$ $\in$ $E(G)$ or $(2)$ $u=u'$ and $vv'$ $\in$ $E(H)$.
\item The \textit{corona} of $G$ with $H$, denoted $G$ $\circ$ $H$, is the graph of order $|G||H|$ $+$ $|G|$ obtained by taking one copy of $G$ and $|G|$ copies of $H$, and joining all the vertices in the \textit{i}th copy of $H$ to the \textit{i}th vertex of $G$.
\end{itemize}

\indent Zero forcing is an iterative process of forcing vertices of a graph $G$ to be blue under a specified color-change rule. This topic has been studied over the last ten years because of its application to different research areas such as combinatorial matrix theory (particularly in minimum rank problems), network algorithms, quantum systems, etc. The zero forcing number of a graph was introduced in \cite{AIM} as an upper bound for the maximum nullity of a graph. In \cite{FJS}, the failed-type parameter of zero forcing was first considered. The notion of the failed zero forcing number of a graph is a relatively new topic and thus poses a lot of open problems. 

\indent This paper is organized as follows. In Section \ref{sec2}, we introduce failed zero forcing sets and the failed zero forcing number of a graph $G$, denoted $F(G)$. We also present some known results from \cite{FJS}. In Section \ref{sec3}, we construct maximal failed zero forcing sets for several graph products. In particular, we consider products of two paths, two cycles, or two complete graphs. We also show that the cardinality of such maximal sets is a sharp bound for $F(G)$ for each type of product. 

\section{Failed zero forcing sets and the graph parameter $F(G)$}\label{sec2}

We begin with a graph $G$ with each vertex of $G$ having an initial coloring of either blue or white. We let $S$ to be the set of all blue vertices. Given a color-change rule, the \textit{derived coloring} is the result of iteratively applying the color-change rule on $G$ with $S$ as the initial set until no more changes are possible. 

\begin{definition}
\textit{Standard zero forcing} 
\begin{itemize}
\item \textit{Standard color-change rule:} {\normalfont If $v \in V(G) \backslash S$ is the unique white vertex of $G$ that is adjacent to a blue vertex $u$ in $S$, then change the color of $v$ from white to blue.}
\item \textit{Iteration:} {\normalfont Update $S$ by adjoining all vertices made blue.}
\item {\normalfont We say that a subset $S$ of $V(G)$ is a \textit{zero forcing set} for $G$ if the derived coloring under the standard color-change rule with $S$ as the initial set is all blue. Otherwise, we say that $S$ is a \textit{failed zero forcing set.}}
\item {\normalfont The \textit{zero forcing number of a graph $G$}, denoted $Z(G)$, is the minimum of $|Z|$ taken over all zero forcing sets for $Z$ for $G$.}
\item {\normalfont The \textit{failed zero forcing number of a graph $G$}, denoted $F(G)$, is the maximum of $|F|$ taken over all failed zero forcing sets $F$ for $G$.}
\end{itemize}
\end{definition}

\begin{observation}\label{obs2.2}
Let $S$ be a zero forcing set for a graph $G$. If $S'$ is a subset of $V(G)$ that contains $S$, then $S'$ is also a zero forcing set for $G$. 
\end{observation}

\begin{observation}\label{obs2.3}
Let $S$ be a failed zero forcing set for a graph $G$. If $S'$ is a subset of $V(G)$ that is contained in $S$, then $S'$ is also a failed zero forcing set for $G$. 
\end{observation}

\indent Observations \ref{obs2.2} and \ref{obs2.3} justify why it is natural and intuitive to consider the minimum size of a (positive semidefinite) zero forcing set and the maximum size of a failed zero forcing set. A zero forcing set with cardinality equal to the zero forcing number is said to be a \textit{minimum zero forcing set}. On the other hand, a failed zero forcing set with cardinality equal to the failed zero forcing number is said to be a \textit{maximum failed zero forcing set}. \\
\indent Our main interest, however, lies on the construction of maximal failed zero forcing sets. We say that a set $F$ is maximal with respect to being failed if $F$ is a failed zero forcing set and any set properly containing $F$ is already a zero forcing set. 

\begin{definition}\label{def2.5} 
Let $G$ be a graph. We say that a proper subset $S$ of $V(G)$ is \textit{stalled} if no color changes are possible from $S$ under the standard color-change rule. 
\end{definition}

\begin{observation}\label{obs2.6}
Maximal and maximum failed zero forcing sets are stalled under the standard color-change rule. Furthermore, a maximum failed zero forcing set is maximal. 
\end{observation}

We next enumerate some existing results from \cite{FJS} on failed zero forcing including a property of a failed zero forcing set and values of $F(G)$ for several graph families. On products of graphs, a lower bound for the failed zero forcing number of a Cartesian product of two graphs is given. 

\begin{theorem}\label{thm2.7}
\cite[Observation $2.1$, Theorem $2.4$]{FJS} For any graph $G$ with $n$ vertices, $Z(G) - 1 \leq F(G) \leq n - 1$. $F(G) = n -1$ if and only if $G$ has an isolated vertex. If in addition, $G$ is connected, then $F(G) = n-2$ if and only if $G$ has a module of order 2. 
\end{theorem}

\begin{theorem}\label{thm2.8}
\cite[Theorems $3.1$, $3.2$, $3.3$, $3.4$, $3.6$, and $3.8$]{FJS} Let $G$ be a simple, finite, and connected graph. Then
$$F(G) =
\begin{cases}
\lceil \frac{n-2}{2} \rceil, & \text{if} \hspace{2mm} G = P_n \\ 
\lfloor \frac{n}{2} \rfloor, & \text{if} \hspace{2mm} G = C_n, n \geq 3 \\
n - 2, & \text{if} \hspace{2mm} G = K_n, n \geq 2 \\
n - 2, & \text{if} \hspace{2mm} G \hspace{1.7mm} \text{is an m-ary tree on n vertices.} \hspace{1.7mm} m \geq 2 \\
3, & \text{if} \hspace{2mm} G = W_n, n = 5\\
\lfloor \frac{2n-2}{3} \rfloor, & \text{if} \hspace{2mm} G = W_n, n \geq 4, n \neq 5 \\
m+n-2 & \text{if} \hspace{2mm} G = K_{m,n}, m+n \geq 3, m \geq n \\
6, & \text{if} \hspace{2mm} G \hspace{1.7mm} \text{is a Petersen graph.}   
\end{cases}$$
\end{theorem}

\begin{theorem}\label{thm2.9}
\cite[Theorem $4.1$]{FJS} Let $G$ and $H$ be graphs with $n$ and $m$ vertices, respectively. Then $F(G$ $\square$ $H)$ $\geq$ max$\{mF(G), nF(H)\}$. 
\end{theorem}

\section{Failed zero forcing number of product graphs}\label{sec3}

\subsection{Cartesian products}
\setcounter{subsection}{0}
The next proposition immediately follows from Theorems \ref{thm2.8} and \ref{thm2.9}.

\begin{proposition}\label{prop3.1}
Let $\mathscr{G}$ be the Cartesian product of $G$ and $H$, where $G$ has $n$ vertices. \\
1. If $H = P_t$, then, $F(\mathscr{G})$ $\geq$ max $\big\{tF(G), n \lceil \frac{t-2}{2} \rceil \big\}$ \\
2. If $H = C_t$, then, $F(\mathscr{G})$ $\geq$ max $\big\{tF(G), n \lfloor \frac{t}{2} \rfloor \big\}$ \\
3. If $H = K_t$, then, $F(\mathscr{G})$ $\geq$ max $\big\{tF(G), n(t-2)\big\}$ 
\end{proposition}

\begin{theorem}\label{thm3.2}
Let $G = P_n$ $\square$ $P_m$, where $n, m \in \mathbb{Z}$ and $n \geq m \geq 2$. A failed zero forcing set for $P_n$ $\square$ $P_m$, which is maximal if $m > 2$, can be constructed using the following algorithm: \\
\normalfont{
\noindent \textit{Input:} $G=$ $P_n$ $\square$ $P_m$, where $n$ and $m$ are the number of vertices of $P_n$ and $P_m$, respectively. \\
\noindent \textit{Output:} A failed zero forcing set $F$ for $P_n$ $\square$ $P_m$, which is maximal if $m > 2$. \\
\noindent \textit{Algorithm:} \\
1. Label each vertex of $G$ as $v_{i,j}$, where $i$ and $j$ denote which copy of $P_n$ and which copy of $P_m$ the vertex is in, respectively. \\
2. Let $W_0 = \{v_{i,i} : i = 1, \ldots m\}$. That is, we start from vertex $v_{1, 1}$ and move diagonally to $v_{m,m}$. \\
3. Set $W = W_0$, $v = v_{m,m}$, and $k=1$. \\
4. While $v$ is in a diagonal that does not contain a vertex in $W$ $\backslash$ $\{v\}$: \\ 
(a) $W_k = \emptyset$. \\
(b) Move away from $v$ along such a diagonal and add every vertex to $W_k$ until a vertex located in the first row, last row, or last column has been added to $W_k$. \\
(c) Let $v'$ be the last vertex added to $W_k$. \\ 
(d) Update $v$ with $v = v'$ and $W$ with $W = W \cup W_k$. \\
(e) Update $k$ with $k=k+1$. \\
5. Consider the two diagonals where $v$ belongs. While there is a vertex not yet in $W$, move away from $v$ along the diagonal different from the diagonal last traced and add every vertex to $W$. Let $W_{k+1}$, $W_{k+2}$, and $W_{k+3}$ be the set (possibly empty) of vertices added to $W$ from each diagonal considered. \\
6. Set $F = V(G)$ $\backslash$ $W$. $F$ is now a failed zero forcing set for $G$ and is maximal if $m > 2$.}
\end{theorem}

\proof
Let $F$ be the set of vertices obtained from the algorithm. Color all vertices in $F$ blue and all vertices outside $F$ white. The process of selecting which vertices are to be excluded from $F$ assures that each blue vertex with a white neighbor will have at least two white neighbors. If we add a white vertex to $F$, the process of selecting also assures that one white vertex will then be a unique white neighbor of a blue vertex, causing a color-change rule. Iteratively, this results to coloring all vertices blue. Hence, $F$ is a maximal failed zero forcing set. 
\qed

\begin{corollary}\label{cor3.3}
Let $G$ $=$ $P_n$ $\square$ $P_m$, where $n, m \in \mathbb{Z}$ with $n \geq m \geq 2$. Let $r$ be the remainder when $n-m$ is divided by $m-1$ unless $n=m$, in which case let $r=0$. Then
$$F(G) \geq \begin{cases} 
nm-n, & \; if \; r = 0 \\
nm-n-m+2, & \; if \; r \neq 0 \\
\end{cases}$$
\end{corollary}

\proof
Let $G$ $=$ $P_n$ $\square$ $P_m$. Construct a maximal failed zero forcing set $F$ for $G$ using Theorem \ref{thm3.2}. Let $k = \lfloor \frac{n-m}{m-1} \rfloor$. Suppose $r=0$. Then $F_{k+1} = F_{k+2} = F_{k+3} = \emptyset$. Hence, $F(G) \geq$ $|F|$ $=$ $|V(G)|$ $-$ $\sum\limits_{t=0}^{k} F_t$ $=$ $nm - (m + \lceil \frac{k}{2} \rceil(m-1)+ \lfloor \frac{k}{2} \rfloor (m-1))$ $=$ $nm - (m + k(m-1))$ $=$ $nm - n$. Suppose $r \neq 0$. Then, $F(G) \geq$ $|F|$ $=$ $|V(G)|$ $-$ $\sum\limits_{t=0}^{k+3} F_t$ $=$ $nm - (m + \lceil \frac{k}{2} \rceil(m-1)+ \lfloor \frac{k}{2} \rfloor (m-1) + r + m-r-1 + r-1)$ $=$ $nm - (m + k(m-1) + m + r - 2)$ $=$ $nm - (m+k(m-1) + r) - m + 2$ $=$ $nm - n - m + 2$.
\qed

\begin{corollary}\label{cor3.4}
\cite[Theorem $4.2$]{FJS} For the square grid graph $P_n$ $\square$ $P_n$ with $n \geq 2$, $F(P_n$ $\square$ $P_n)$ $=$ $n^2 -n$. 
\end{corollary}
\proof
It immediately follows from Corollary \ref{cor3.3} that $F(P_n$ $\square$ $P_m)$ $\geq$ $n^2 -n$. Equality is proved using the pigeonhole principle and by contradiction when a set containing more than $n^2-n$ vertices is assumed to be failed. 
\qed

\begin{theorem}\label{thm3.5} 
Let $m$ and $n$ be integers greater than 2. If $G = C_m$ $\square$ $C_n$, then $F(G) \geq \lceil \frac{mn}{2} \rceil$.
\end{theorem}
\proof
Let $G = C_m$ $\square$ $C_n$ and $H_i$ be the induced subgraph of the $i$th copy of $C_n$ for $i = 1, ..., m$. Label each vertex of $G$ as $v_{i,j}$ where $i$ and $j$ denote the copy of $C_n$ and the copy of $C_m$ where the vertex is in, respectively. Denote $v_{0,j}$ as $v_{m,j}$, $v_{m+1,j}$ as $v_{1,j}$, $v_{i,0}$ as $v_{i,n}$, and $v_{i,n+1}$ as $v_{i,1}$. Let $S$ be a set of vertices of $G$ such that $v_{i,j}$ is in $S$ if and only if $i$ and $j$ have the same parity. Consider $v_{i,j} \in S$. Then $v_{i,j}$ is adjacent to $v_{i-1,j}$ and $v_{i+1,j}$. Since $i$ and $j$ have the same parity, $j$ has a parity opposite to that of $i-1$ and $i+1$. Hence, $v_{i-1,j}$ and $v_{i+1,j}$ are both not in $S$. This implies that every vertex in $S$ has at least two neighbors outside in $S$. Hence, $S$ is stalled under the standard-color-change rule and $F(G) \geq S$. \\
\indent Let $i \in \{1, ..., m\}$. If $n$ is even, there are $\frac{n}{2}$ values of $j$ in $\{1, ..., n\}$ for which $i$ and $j$ have the same parity. Hence, $|S| = \frac{mn}{2}$. On the other hand, let $j \in \{1, ..., n\}$. If $m$ is even, there are $\frac{m}{2}$ values of $i$ in $\{1, ..., n\}$ for which $i$ and $j$ have the same parity. Hence, $|S| = \frac{nm}{2}$. Thus, if at least one of $m$ and $n$ are even, $F(G) \geq \frac{mn}{2} = \lceil \frac{mn}{2} \rceil$. \\
\indent Suppose both $n$ and $m$ are odd. Then the number of odd $j$'s in $\{1,...,n\}$ is greater than the number of even $j$'s in $\{1,...,n\}$. Similarly, the number of odd $i$'s in $\{1,...,m\}$ is greater than the number of even $i$'s in $\{1,...,m\}$. This implies that there are $\lceil \frac{n}{2} \rceil$ odd $j$'s, $\lceil \frac{n}{2} \rceil$ even $j$'s, $\lceil \frac{m}{2} \rceil$ odd $i$'s, and $\lceil \frac{n}{2} \rceil$ even $i$'s. Consider $i \in \{1, ..., m\}$. If $i$ is odd, then there are $\lceil \frac{n}{2} \rceil$ values of $j$ for which $v_{i,j} \in S$. Hence, there are $\lceil \frac{m}{2} \rceil$ $\lceil \frac{n}{2} \rceil$ vertices $v_{i,j}$ in $S$ for which $i$ and $j$ are both odd. If $i$ is even, then there are $\lfloor \frac{n}{2} \rfloor$ values of $j$ for which $v_{i,j} \in S$. Hence, there are $\lfloor \frac{m}{2} \rfloor$ $\lfloor \frac{n}{2} \rfloor$ vertices $v_{i,j}$ in $S$ for which $i$ and $j$ are both even. Thus, $|S| =   \lceil \frac{m}{2} \rceil \lceil \frac{n}{2} \rceil + \lfloor \frac{m}{2} \rfloor \lfloor \frac{n}{2} \rfloor$. Since $n$ and $m$ are both odd, $|S| = \big(\frac{m+1}{2}\big)\big(\frac{n+1}{2}\big) + \big(\frac{m-1}{2}\big)\big(\frac{n-1}{2}\big) = \frac{mn+1}{2} = \lceil \frac{mn}{2} \rceil$. 
\qed

\begin{theorem}\label{thm3.6} 
Let $G = P_2$ $\square$ $C_n$, where $n$ is an integer greater than 2. Then, $F(G) \geq \lceil \frac{n}{2} \rceil + 3 \lfloor \frac{n}{4} \rfloor$.  
\end{theorem}
\proof Let $G= P_2$ $\square$ $C_n$. Label each vertex in $G$ as $v_{i,j}$, where $i$ denotes which copy of $C_n$ and $j$ denotes which copy of $P_2$ the vertex is in. Denote $v_{0,j}$ as $v_{2,j}$, $v_{3,j}$ as $v_{1,j}$, $v_{i,0}$ as $v_{i,n}$, and $v_{i,n+1}$ as $v_{i,1}$. \\
\indent Suppose $n = 4k+l$, where $k \in \mathbb{N} \cup \{0\}$ and $l \in \{0,1, 2, 3\}$. \\
\indent Case 1: $l = 0$. \\
\indent Construct $S$ by adding to $S$ the following vertices: \\
\indent 1. vertices $v_{1,j_1}$, where $j_1$ is not divisible by 4.  \\
\indent 2. vertices $v_{2,j_2}$, where $j_2$ is even. \\
\indent Color all vertices in $S$ blue. We claim that $S$ is failed. Indeed, suppose $v \in S$. Then $v$ is of the form $v_{1, 4k_1 + 1}$, $v_{1,4k_1 + 2}$, $v_{1,4k_1 + 3}$, or $v_{2,2k_2}$, where $k_1$ and $k_2$ are integers such that $0 \leq k_1 < k$ and $0 \leq k_2 \leq 2k$. Now, $v_{1, 4k_1 + 1}$ is adjacent to $v_{2,4k_1+1}$ and $v_{1,4k_1}$, both of which are white. $v_{1,4k_1+2}$ has no white neighbors. $v_{1, 4k_1 + 3}$ is adjacent to $v_{2,4k_1+3}$ and $v_{1,4(k_1+1)}$, both of which are white. Lastly, $v_{2,2k_2}$ is adjacent to $v_{2,2k_2-1}$ and $v_{2,2k_2+1}$, which are both white. Therefore, each vertex in $S$ has either zero or two white neighbors. Hence, $S$ is stalled and failed. \\
\indent If Case 1 happens, $F(G) \geq |S| = 3 (\frac{n}{4}) + \frac{n}{2}$ $=$ $3 \lfloor \frac{n}{4} \rfloor$ $+$ $\lceil \frac{n}{2} \rceil $. \\
\indent Case 2: $l = 1$. \\
\indent Construct $S$ by adding to $S$ the following vertices: \\
\indent 1. vertices $v_{1,j_1}$, where $j_1$ is not divisible by 4 and $j_1 \neq n$.  \\
\indent 2. vertices $v_{2,j_2}$, where $j_2$ is even. \\
\indent 3. vertex $v_{2,n}$ \\
\indent Color all vertices in $S$ blue. We claim that $S$ is failed. Indeed, suppose $v \in S$. Then $v$ is of the form $v_{1, 4k_1 + 1}$, $v_{1,4k_1 + 2}$, $v_{1,4k_1 + 3}$, $v_{2,2k_2}$ or $v_{2,n}$, where $k_1$ and $k_2$ are integers such that $0 \leq k_1 < k$ and $0 \leq k_2 \leq 2k$. Now, $v_{1, 4k_1 + 1}$ is adjacent to $v_{2,4k_1+1}$ and $v_{1,4k_1}$, both of which are white. The vertex $v_{1,4k_1+2}$ has no white neighbors. The vertex $v_{1, 4k_1 + 3}$ is adjacent to $v_{2,4k_1+3}$ and $v_{1,4(k_1+1)}$, both of which are white. On the other hand, $v_{2,2k_2}$ is adjacent to $v_{2,2k_2-1}$ and $v_{2,2k_2+1}$, which are both white. Lastly, $v_{2,n}$ is adjacent to the white vertices $v_{1,n}$ and $v_{2,1}$. Therefore, each vertex in $S$ has either zero or two white neighbors. This implies that $S$ is stalled and failed. \\
\indent If Case 2 happens, $F(G) \geq |S| = 3 \lfloor \frac{n}{4} \rfloor + \frac{n-1}{2} + 1$ $=$ $3 \lfloor \frac{n}{4} \rfloor$ $+$ $\frac{n+1}{2}$ $=$ $3 \lfloor \frac{n}{4} \rfloor$ $+$ $\lceil \frac{n}{2} \rceil $. \\
\indent Case 3: $l = 2$. \\
\indent Construct $S$ by adding to $S$ the following vertices: \\
\indent 1. vertices $v_{1,j_1}$, where $j_1$ is not divisible by 4 and $j_1 \neq n, n-1$.  \\
\indent 2. vertices $v_{2,j_2}$, where $j_2$ is even. \\
\indent Color all vertices in $S$ blue. We claim that $S$ is failed. Indeed, suppose $v \in S$. Then $v$ is of the form $v_{1, 4k_1 + 1}$, $v_{1,4k_1 + 2}$, $v_{1,4k_1 + 3}$, or $v_{2,2k_2}$. where $k_1$ and $k_2$ are integers such that $0 \leq k_1 < k$ and $0 \leq k_2 \leq 2k+1$. Now, $v_{1, 4k_1 + 1}$ is adjacent to $v_{2,4k_1+1}$ and $v_{1,4k_1}$, both of which are white. The vertex $v_{1,4k_1+2}$ has no white neighbors. The vertex $v_{1, 4k_1 + 3}$ is adjacent to $v_{2,4k_1+3}$ and $v_{1,4(k_1+1)}$, both of which are white. Lastly, $v_{2,2k_2}$ is adjacent to $v_{2,2k_2-1}$ and $v_{2,2k_2+1}$, which are both white. Therefore, each vertex in $S$ has either zero or two white neighbors. Hence, $S$ is stalled and failed. \\
\indent If Case 3 happens, $F(G) \geq |S|$ $=$ $3(\frac{n-2}{4}) + \frac{n}{2}$ $=$ $3 \lfloor \frac{n}{4} \rfloor$ $+$ $\lceil \frac{n}{2} \rceil $. \\
\indent Case 4: $l = 3$. \\
\indent Construct $S$ by adding to $S$ the following vertices: \\
\indent 1. vertices $v_{1,j_1}$, where $j_1$ is not divisible by 4 and $j_1 \neq n, n-1$.  \\
\indent 2. vertices $v_{2,j_2}$, where $j_2$ is even. \\
\indent Color all vertices in $S$ blue. We claim that $S$ is failed. Indeed, suppose $v \in S$. Then $v$ is of the form $v_{1, 4k_1 + 1}$, $v_{1,4k_1 + 2}$, $v_{1,4k_1 + 3}$, $v_{2,2k_2}$, or $v_{1,n-2}$, where $k_1$ and $k_2$ are integers such that $0 \leq k_1 < k$ and $0 \leq k_2 \leq 2k+1$. The vertex $v_{1, 4k_1 + 1}$ is adjacent to $v_{2,4k_1+1}$ and $v_{1,4k_1}$, both of which are white. The vertex $v_{1,4k_1+2}$ has no white neighbors. The vertex $v_{1, 4k_1 + 3}$ is adjacent to $v_{2,4k_1+3}$ and $v_{1,4(k_1+1)}$, both of which are white. On the other hand, $v_{2,2k_2}$ is adjacent to $v_{2,2k_2-1}$ and $v_{2,2k_2+1}$, which are both white. Lastly, $v_{1,n-2}$ is adjacent to the white vertices $v_{1,n-1}$, $v_{1,n-3}$, and $v_{2,n-2}$. Therefore, each vertex in $S$ has either zero or at least two white neighbors. This implies that $S$ is stalled and failed. \\
\indent If Case 4 happens, $F(G) \geq |S| = 3 \lfloor \frac{n}{4} \rfloor + 1 + \frac{n-1}{2}$ $=$ $3 \lfloor \frac{n}{4} \rfloor$ $+$ $\frac{n+1}{2}$ $=$ $3 \lfloor \frac{n}{4} \rfloor$ $+$ $\lceil \frac{n}{2} \rceil $. 
\qed

\begin{theorem}\label{thm3.7} 
\cite[Theorem $4.3$]{FJS} Let $n, m \in \mathbb{N}$ such that $n \geq 4$ and $m \geq 2$. Then, $F(K_n$ $\square$ $K_m)$ $= nm-4$. 
\end{theorem}

\setcounter{subsection}{1}
\subsection{Strong products}
\setcounter{subsection}{7}

\begin{theorem}\label{thm3.8}
Let $G$ and $H$ be graphs with $n$ and $m$ vertices, respectively. Let \linebreak $S \subseteq V(G$ $\boxtimes$ $H)$. If every vertex in $S$ is adjacent to at least two vertices in \linebreak $V(G$ $\square$ $H)$ $\backslash$ $S$, then $F(G$ $\boxtimes$ $H)$ $\geq |S|$. In particular, if the cardinality of such set $S$ is $F(G$ $\square$ $H)$, then $F(G$ $\boxtimes$ $H)$ $\geq$ max$\{mF(G), nF(H)\}$.
\end{theorem}
\proof
Let $G$ and $H$ be graphs with $n$ and $m$ vertices, respectively. Let $S$ be a set of vertices of the graph $G$ $\square$ $H$. To construct $G$ $\boxtimes$ $H$, we simply add more edges to $G$ $\square$ $H$. Doing this retains or increases the number of neighbors a vertex in $S$ has. Suppose all vertices in $S$ have at least two white neighbors in $V(G$ $\square$ $H)$ $\backslash$ $S$. Then, all vertices in $S$ will also have at least two white neighbors in $V(G$ $\boxtimes$ $H)$ $\backslash$ $S$. Thus, no standard color changes are possible from $S$. Therefore, $S$ is stalled and $S$ is a failed zero forcing set for $G$ $\boxtimes$ $H$. Hence, $F(G$ $\boxtimes$ $H)$ $\geq |S|$. \\
\indent Suppose a maximum failed zero forcing set $S$ for $G$ $\square$ $H$ is a set whose every element has at least two white neighbors outside the set in $G$ $\square$ $H$. By the preceding paragraph, $F(G$ $\boxtimes$ $H)$ $\geq$ $|S|$ $=$ $F(G$ $\square$ $H)$. Finally, by Theorem \ref{thm2.9}, \linebreak $F(G$ $\boxtimes$ $H)$ $\geq$ max$\{mF(G), nF(H)\}$
\qed

\begin{lemma}\label{lemma3.9}
Let $G$ be the path $P_n$ with $n$ vertices. If $V(G) = \{v_1, ..., v_n\}$ and $E(G) = \{v_iv_{i+1} : i = 1, ..., n-1\}$, then the following sets are standard zero forcing sets for $G$: 
\begin{enumerate}
\item $S = \{v\}$, where $v$ is an endpoint of $P_n$, that is, $v = v_1$ or $v = v_n$
\item $S = \{v_i, v_j\}$, where $v_i$ and $v_j$ are consecutive vertices of $P_n$, that is, $|i-j|=1$ for some $i, j \in \{1, ..., n\}$. 
\end{enumerate}
\end{lemma}
\proof
Suppose $S = \{v_1\}$. Applying the standard color-change rule, we have the following chain of forces:
$$v_1 \rightarrow v_2 \rightarrow \cdots \rightarrow v_n$$
\indent On the other hand, if $S = \{v_n\}$, we have the following chain of forces:
$$v_n \rightarrow v_{n-1} \rightarrow \cdots \rightarrow v_2$$
\indent Hence, a set of an endpoint of $P_n$ is a standard zero forcing set for $P_n$. Now, suppose $S = \{v_i, v_{i+1}\}$ for some $i \in \{1, ..., n\}$. Since $\{v_1\}$ and $\{v_n\}$ are standard zero forcing sets for $P_n$, by Observation \ref{obs2.2}, $S$ is also a standard zero forcing set for $P_n$ when $i=1$ or $i=n-1$. If $i \neq 1$ and $i \neq n-1$, the following are the chain of forces:
$$v_i \rightarrow v_{i-1} \rightarrow \cdots v_1$$ 
$$v_{i+1} \rightarrow v_{i+2} \rightarrow \cdots v_n$$
\indent Therefore, a set of two consecutive vertices of $P_n$ is a standard zero forcing set for $P_n$.
\qed

\begin{theorem}\label{theorem3.10} 
If $n$ and $m$ are two integers such that $n \geq m \geq 2$, then \linebreak $F(P_n$ $\boxtimes$ $P_m)$ $\geq$ $s$ and there exists a set $S$ with cardinality $s$ that is maximal and failed, where $s = nm - m + \lceil \frac{m-4}{3} \rceil$. 
\end{theorem}
\proof
Let $v_{i,j}$ be the vertex of $F(P_n$ $\boxtimes$ $P_m)$ corresponding to the $i$th copy of $P_m$ and $j$th copy of $P_n$. If $m=2$, then $G = P_n$ $\boxtimes$ $P_2$. Examining the graph $G$ below, $v_{1,1}$ and $v_{1,2}$ form a module of order 2 since both are adjacent to only each other and the vertices $v_{2,1}$ and $v_{2,2}$. Since $G$ is connected, by Theorem \ref{thm2.7}, $F(G) = |G| - 2 = 2n-2$ $=$ $2n - 2 + \lceil \frac{2-4}{3} \rceil$. Hence, our bound is sharp when $m=2$. Thus, a failed set with cardinality $2n-2$ is maximum and therefore, maximal.
\begin{figure}[H] 
\begin{center}

\begin{tikzpicture}[line cap=round,line join=round,scale=0.56]
\draw [line width=1pt] (-4,4)-- (-4,2);
\draw [line width=1pt] (-2,4)-- (-4,2);
\draw [line width=1pt] (-4,2)-- (-2,2);
\draw [line width=1pt] (-2,2)-- (-2,4);
\draw [line width=1pt] (-2,4)-- (-4,4);
\draw [line width=1pt] (-4,4)-- (-2,2);
\draw [line width=0.5pt,dash pattern=on 2pt off 1pt] (-2,4)-- (0,4);
\draw [line width=0.5pt,dash pattern=on 2pt off 1pt] (-2,2)-- (0,2);
\draw [line width=0.5pt,dash pattern=on 2pt off 1pt] (0,4)-- (2,4);
\draw [line width=0.5pt,dash pattern=on 2pt off 1pt] (0,2)-- (2,2);
\draw [line width=0.5pt,dash pattern=on 1pt off 1pt] (-2,4)-- (0,2);
\draw [line width=0.5pt,dash pattern=on 1pt off 1pt] (2,4)-- (0,2);
\draw [line width=0.5pt,dash pattern=on 1pt off 1pt] (0,4)-- (2,2);
\draw [line width=0.5pt,dash pattern=on 1pt off 1pt] (-2,2)-- (0,4);
\draw [line width=1pt] (2,4)-- (2,2);
\draw (-4.28,4.64) node[anchor=north west] {$v_{11}$};
\draw (-2.3,4.64) node[anchor=north west] {$v_{21}$};
\draw (1.7,4.64) node[anchor=north west] {$v_{n,1}$};
\draw (1.7,1.86) node[anchor=north west] {$v_{n,2}$};
\draw (-4.28,1.91) node[anchor=north west] {$v_{12}$};
\draw (-2.3,1.91) node[anchor=north west] {$v_{22}$};
\draw [fill=white] (-4,4) circle (0.2);
\draw [fill=white] (-4,2) circle (0.2);
\draw [fill=blue] (-2,4) circle (0.2);
\draw [fill=blue] (-2,2) circle (0.2);
\draw [fill=blue] (2,4) circle (0.2);
\draw [fill=blue] (2,2) circle (0.2);
\end{tikzpicture}
\end{center}
\caption{\label{fig3.16} $P_n$ $\boxtimes$ $P_2$.}
\end{figure}
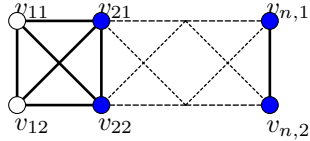
If $m=3$, then $G = P_n$ $\boxtimes$ $P_3$. Examining the graph $G$ below, $G$ has no module of order 2 nor an isolated vertex. Since $G$ is connected, by Theorem \ref{thm2.7}, $F(G) < |G| - 2 = 3n - 2$. The set of $3n-3$ blue vertices, as seen in the graph below, is stalled and thus failed. Hence, $F(G) = 3n-3$ $=$ $3n - 2 + \lceil \frac{3-4}{3} \rceil$. Hence, our bound is sharp when $m=3$. Thus, a failed set with cardinality $3n-3$ is maximum and therefore, maximal. 
\begin{figure}[H] 
\begin{center}
\begin{tikzpicture}[line cap=round,line join=round,scale=0.56]
\draw [line width=1pt] (-4,4)-- (-4,2);
\draw [line width=1pt] (-4,2)-- (-2,2);
\draw [line width=1pt] (-2,2)-- (-2,4);
\draw [line width=1pt] (-2,4)-- (-4,4);
\draw [line width=1pt] (-4,4)-- (-2,2);
\draw [line width=0.5pt,dash pattern=on 2pt off 1pt] (-2,4)-- (0,4);
\draw [line width=0.5pt,dash pattern=on 2pt off 1pt] (-2,2)-- (0,2);
\draw [line width=0.5pt,dash pattern=on 2pt off 1pt] (0,4)-- (2,4);
\draw [line width=0.5pt,dash pattern=on 2pt off 1pt] (0,2)-- (2,2);
\draw [line width=1pt] (-2,4)-- (-4,2);
\draw [line width=1pt,dash pattern=on 1pt off 1pt] (-2,4)-- (0,2);
\draw [line width=1pt,dash pattern=on 1pt off 1pt] (2,4)-- (0,2);
\draw [line width=1pt,dash pattern=on 1pt off 1pt] (0,4)-- (2,2);
\draw [line width=1pt,dash pattern=on 1pt off 1pt] (-2,2)-- (0,4);
\draw [line width=1pt] (2,4)-- (2,2);
\draw [line width=1pt] (-4,0)-- (-2,0);
\draw [line width=0.5pt,dash pattern=on 2pt off 1pt] (-2,0)-- (0,0);
\draw [line width=0.5pt,dash pattern=on 2pt off 1pt] (0,0)-- (2,0);
\draw [line width=1pt] (-4,2)-- (-4,0);
\draw [line width=1pt] (-2,2)-- (-2,0);
\draw [line width=1pt] (-4,2)-- (-2,0);
\draw [line width=1pt] (-2,2)-- (-4,0);
\draw [line width=1pt] (2,2)-- (2,0);
\draw [line width=1pt,dash pattern=on 1pt off 1pt] (-2,2)-- (0,0);
\draw [line width=1pt,dash pattern=on 1pt off 1pt] (0,2)-- (-2,0);
\draw [line width=1pt,dash pattern=on 1pt off 1pt] (0,2)-- (2,0);
\draw [line width=1pt,dash pattern=on 1pt off 1pt] (2,2)-- (0,0);
\draw [fill=white] (-4,4) circle (0.2);
\draw [fill=white] (-4,2) circle (0.2);
\draw [fill=blue] (-2,4) circle (0.2);
\draw [fill=blue] (-2,2) circle (0.2);
\draw [fill=blue] (2,4) circle (0.2);
\draw [fill=blue] (2,2) circle (0.2);
\draw [fill=white] (-4,0) circle (0.2);
\draw [fill=blue] (-2,0) circle (0.2);
\draw [fill=blue] (2,0) circle (0.2);
\end{tikzpicture}
\end{center}
\caption{\label{fig3.17} $P_n$ $\boxtimes$ $P_3$.}
\end{figure}
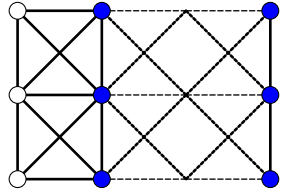
\indent Suppose $m > 3$. Let $G = P_n$ $\boxtimes$ $P_m$. Construct $S$ as follows by adding to $S$ all vertices $v_{i,j}$ such that $2 \leq i \leq n$ and $1 \leq j \leq m$. If, in addition, $m > 4$, add to $S$ the vertices $v_{1,3k}$ for $k=1, ..., \lceil \frac{m-4}{3} \rceil$. \\ 
\indent We claim that $S$ is failed. Indeed, suppose $v_{i,j} \in S$. \\ 
\indent Case 1: Suppose $v_{1,j} \in S$. Then $v_{1,j}$ is adjacent to $v_{1,j-1}$ and $v_{1,j+1}$, which are both white. 
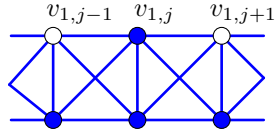
\begin{figure}[H]
\begin{center}
\begin{tikzpicture}[line cap=round,line join=round,scale=0.56]
\draw [line width=1pt,color=blue] (-3,4)-- (-3,2);
\draw [line width=1pt,color=blue] (-3,4)-- (-1,2);
\draw [line width=1pt,color=blue] (-1,2)-- (-1,4);
\draw [line width=1pt,color=blue] (-1,4)-- (1,2);
\draw [line width=1pt,color=blue] (1,2)-- (1,4);
\draw [line width=1pt,color=blue] (1,4)-- (-1,2);
\draw [line width=1pt,color=blue] (-3,2)-- (-1,4);
\draw [line width=1pt,color=blue] (-3,4)-- (-1,4);
\draw [line width=1pt,color=blue] (-3,2)-- (-1,2);
\draw [line width=1pt,color=blue] (-1,4)-- (1,4);
\draw [line width=1pt,color=blue] (-1,2)-- (1,2);
\draw (-1.3,4.96) node[anchor=north west] {$v_{1,j}$};
\draw (-3.38,4.96) node[anchor=north west] {$v_{1,j-1}$};
\draw (0.5,4.96) node[anchor=north west] {$v_{1,j+1}$};
\draw [line width=1pt,color=blue] (1,4)-- (2,4);
\draw [line width=1pt,color=blue] (1,2)-- (2,2);
\draw [line width=1pt,color=blue] (-3,4)-- (-4,4);
\draw [line width=1pt,color=blue] (-3,2)-- (-4,2);
\draw [line width=1pt,color=blue] (1,4)-- (2,3);
\draw [line width=1pt,color=blue] (1,2)-- (2,3);
\draw [line width=1pt,color=blue] (-3,4)-- (-4.04,2.84);
\draw [line width=1pt,color=blue] (-3,2)-- (-4.04,2.84);
\draw [fill=white] (-3,4) circle (0.2);
\draw [fill=blue] (-1,4) circle (0.2);
\draw [fill=white] (1,4) circle (0.2);
\draw [fill=blue] (-3,2) circle (0.2);
\draw [fill=blue] (-1,2) circle (0.2);
\draw [fill=blue] (1,2) circle (0.2);
\end{tikzpicture}
\end{center}
\caption{\label{fig3.18} Case 1: $v_{1,j} \in S$}
\end{figure}
\indent Case 2: Suppose $v_{i,j} \in S$ with $i \neq 1$. \\
\indent If $i > 2$, $v_{i,j}$ will only have blue neighbors. Suppose $i=2$. \\
\indent Case 2.1: Consider $v_{2,j}$ for $j=1, 2$. Then $v_{2,j}$ is adjacent to $v_{1,1}$ and $v_{1,2}$, which are both white. 
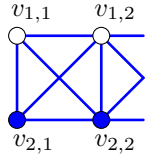
\begin{figure}[H]
\begin{center}
\begin{tikzpicture}[line cap=round,line join=round,scale=0.56]
\draw [line width=1pt,color=blue] (-3,4)-- (-3,2);
\draw [line width=1pt,color=blue] (-3,4)-- (-1,2);
\draw [line width=1pt,color=blue] (-1,2)-- (-1,4);
\draw [line width=1pt,color=blue] (-3,2)-- (-1,4);
\draw [line width=1pt,color=blue] (-3,4)-- (-1,4);
\draw [line width=1pt,color=blue] (-3,2)-- (-1,2);
\draw (-3.32,1.88) node[anchor=north west] {$v_{2,1}$};
\draw (-1.38,4.97) node[anchor=north west] {$v_{1,2}$};
\draw (-3.36,4.97) node[anchor=north west] {$v_{1,1}$};
\draw (-1.38,1.88) node[anchor=north west] {$v_{2,2}$};
\draw [line width=1pt,color=blue] (-1,4)-- (0,4);
\draw [line width=1pt,color=blue] (-1,4)-- (0,3);
\draw [line width=1pt,color=blue] (-1,2)-- (0,3);
\draw [line width=1pt,color=blue] (-1,2)-- (0,2);
\draw [fill=white] (-3,4) circle (0.2);
\draw [fill=white] (-1,4) circle (0.2);
\draw [fill=blue] (-3,2) circle (0.2);
\draw [fill=blue] (-1,2) circle (0.2);
\end{tikzpicture}
\end{center}
\caption{\label{fig3.19} Case $2.1$: $v_{2,1}$ and $v_{2,2}$}
\end{figure}
\indent Case 2.2: Consider $v_{2,j}$ for $j= m-1, m$. Then $v_{2,j}$ is adjacent to $v_{1,m-1}$ and $v_{1,m}$, which are both white. 
\begin{figure}[H]
\begin{center}
\begin{tikzpicture}[line cap=round,line join=round,scale=0.56]
\draw [line width=1pt,color=blue] (-3,4)-- (-3,2);
\draw [line width=1pt,color=blue] (-3,4)-- (-1,2);
\draw [line width=1pt,color=blue] (-1,2)-- (-1,4);
\draw [line width=1pt,color=blue] (-3,2)-- (-1,4);
\draw [line width=1pt,color=blue] (-3,4)-- (-1,4);
\draw [line width=1pt,color=blue] (-3,2)-- (-1,2);
\draw (-3.56,1.88) node[anchor=north west] {$v_{2,m-1}$};
\draw (-1.38,4.97) node[anchor=north west] {$v_{1,m}$};
\draw (-3.56,4.97) node[anchor=north west] {$v_{1,m-1}$};
\draw (-1.38,1.88) node[anchor=north west] {$v_{2,m}$};
\draw [line width=1pt,color=blue] (-3,4)-- (-4,4);
\draw [line width=1pt,color=blue] (-3,4)-- (-4,3);
\draw [line width=1pt,color=blue] (-3,2)-- (-4,3);
\draw [line width=1pt,color=blue] (-3,2)-- (-4,2);
\draw [fill=white] (-3,4) circle (0.2);
\draw [fill=white] (-1,4) circle (0.2);
\draw [fill=blue] (-3,2) circle (0.2);
\draw [fill=blue] (-1,2) circle (0.2);
\end{tikzpicture}
\end{center}
\caption{\label{fig3.20} Case $2.2$: $v_{2,m-1}$ and $v_{2,m}$}
\end{figure}
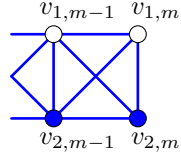
\indent Case 2.3: Suppose $2 < j < m-1$ with $j=3k$ for some $k$ $\in$ $\{1, ..., \lceil \frac{m-4}{3} \rceil \}$. Then, $v_{2,j}$ is adjacent to $v_{1,j-1}$ and $v_{1,j+1}$, which are both white. 
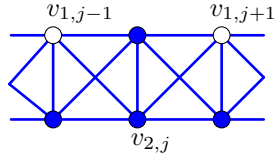
\begin{figure}[H]
\begin{center}
\begin{tikzpicture}[line cap=round,line join=round,scale=0.56]
\draw [line width=1pt,color=blue] (-3,4)-- (-3,2);
\draw [line width=1pt,color=blue] (-3,4)-- (-1,2);
\draw [line width=1pt,color=blue] (-1,2)-- (-1,4);
\draw [line width=1pt,color=blue] (-1,4)-- (1,2);
\draw [line width=1pt,color=blue] (1,2)-- (1,4);
\draw [line width=1pt,color=blue] (1,4)-- (-1,2);
\draw [line width=1pt,color=blue] (-3,2)-- (-1,4);
\draw [line width=1pt,color=blue] (-3,4)-- (-1,4);
\draw [line width=1pt,color=blue] (-3,2)-- (-1,2);
\draw [line width=1pt,color=blue] (-1,4)-- (1,4);
\draw [line width=1pt,color=blue] (-1,2)-- (1,2);
\draw (-1.38,1.88) node[anchor=north west] {$v_{2,j}$};
\draw (-3.38,4.96) node[anchor=north west] {$v_{1,j-1}$};
\draw (0.5,4.96) node[anchor=north west] {$v_{1,j+1}$};
\draw [line width=1pt,color=blue] (1,4)-- (2,4);
\draw [line width=1pt,color=blue] (1,2)-- (2,2);
\draw [line width=1pt,color=blue] (-3,4)-- (-4,4);
\draw [line width=1pt,color=blue] (-3,2)-- (-4,2);
\draw [line width=1pt,color=blue] (1,4)-- (2,3);
\draw [line width=1pt,color=blue] (1,2)-- (2,3);
\draw [line width=1pt,color=blue] (-3,4)-- (-4.04,2.84);
\draw [line width=1pt,color=blue] (-3,2)-- (-4.04,2.84);
\draw [fill=white] (-3,4) circle (0.2);
\draw [fill=blue] (-1,4) circle (0.2);
\draw [fill=white] (1,4) circle (0.2);
\draw [fill=blue] (-3,2) circle (0.2);
\draw [fill=blue] (-1,2) circle (0.2);
\draw [fill=blue] (1,2) circle (0.2);
\end{tikzpicture}
\end{center}
\caption{\label{fig3.21} Case $2.3$: $v_{i,j} \in S$, where $i \neq 1$, $2 < j < m-1$, and $j \equiv 0$ mod $3$}
\end{figure}
\indent Case 2.4: Suppose $2 < j < m-1$ with $j \equiv 1$ mod $3$. Thus, $j+1 \equiv 2$ mod $3$. This implies that $v_{1,j}$ and $v_{1,j+1}$ are both white vertices adjacent to $v_{2,j}$. 
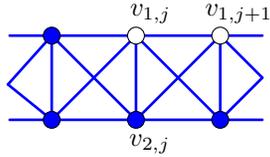
\begin{figure}[H]
\begin{center}
\begin{tikzpicture}[line cap=round,line join=round,scale=0.56]
\draw [line width=1pt,color=blue] (-3,4)-- (-3,2);
\draw [line width=1pt,color=blue] (-3,4)-- (-1,2);
\draw [line width=1pt,color=blue] (-1,2)-- (-1,4);
\draw [line width=1pt,color=blue] (-1,4)-- (1,2);
\draw [line width=1pt,color=blue] (1,2)-- (1,4);
\draw [line width=1pt,color=blue] (1,4)-- (-1,2);
\draw [line width=1pt,color=blue] (-3,2)-- (-1,4);
\draw [line width=1pt,color=blue] (-3,4)-- (-1,4);
\draw [line width=1pt,color=blue] (-3,2)-- (-1,2);
\draw [line width=1pt,color=blue] (-1,4)-- (1,4);
\draw [line width=1pt,color=blue] (-1,2)-- (1,2);
\draw (-1.38,1.88) node[anchor=north west] {$v_{2,j}$};
\draw (0.42,4.97) node[anchor=north west] {$v_{1,j+1}$};
\draw (-1.36,4.97) node[anchor=north west] {$v_{1,j}$};
\draw [line width=1pt,color=blue] (1,4)-- (2,4);
\draw [line width=1pt,color=blue] (1,2)-- (2,2);
\draw [line width=1pt,color=blue] (-3,4)-- (-4,4);
\draw [line width=1pt,color=blue] (-3,2)-- (-4,2);
\draw [line width=1pt,color=blue] (1,4)-- (2,3);
\draw [line width=1pt,color=blue] (1,2)-- (2,3);
\draw [line width=1pt,color=blue] (-3,4)-- (-4.04,2.84);
\draw [line width=1pt,color=blue] (-3,2)-- (-4.04,2.84);
\draw [fill=blue] (-3,4) circle (0.2);
\draw [fill=white] (-1,4) circle (0.2);
\draw [fill=white] (1,4) circle (0.2);
\draw [fill=blue] (-3,2) circle (0.2);
\draw [fill=blue] (-1,2) circle (0.2);
\draw [fill=blue] (1,2) circle (0.2);
\end{tikzpicture}
\end{center}
\caption{\label{fig3.22} Case $2.4$: $v_{i,j} \in S$, where $i \neq 1$, $2 < j < m-1$, and $j \equiv 1$ mod $3$}
\end{figure}
\indent Case 2.5: Suppose $2 < j < m-1$ with $j \equiv 2$ mod $3$. Thus, $j-1 \equiv 1$ mod $3$. This implies that $v_{1,j}$ and $v_{1,j-1}$ are both white vertices adjacent to $v_{2,j}$. 
\begin{figure}[H]
\begin{center}
\begin{tikzpicture}[line cap=round,line join=round,scale=0.56]
\draw [line width=1pt,color=blue] (-3,4)-- (-3,2);
\draw [line width=1pt,color=blue] (-3,4)-- (-1,2);
\draw [line width=1pt,color=blue] (-1,2)-- (-1,4);
\draw [line width=1pt,color=blue] (-1,4)-- (1,2);
\draw [line width=1pt,color=blue] (1,2)-- (1,4);
\draw [line width=1pt,color=blue] (1,4)-- (-1,2);
\draw [line width=1pt,color=blue] (-3,2)-- (-1,4);
\draw [line width=1pt,color=blue] (-3,4)-- (-1,4);
\draw [line width=1pt,color=blue] (-3,2)-- (-1,2);
\draw [line width=1pt,color=blue] (-1,4)-- (1,4);
\draw [line width=1pt,color=blue] (-1,2)-- (1,2);
\draw (-1.38,1.88) node[anchor=north west] {$v_{2,j}$};
\draw (-3.38,4.96) node[anchor=north west] {$v_{1,j-1}$};
\draw (-1.36,4.96) node[anchor=north west] {$v_{1,j}$};
\draw [line width=1pt,color=blue] (1,4)-- (2,4);
\draw [line width=1pt,color=blue] (1,2)-- (2,2);
\draw [line width=1pt,color=blue] (-3,4)-- (-4,4);
\draw [line width=1pt,color=blue] (-3,2)-- (-4,2);
\draw [line width=1pt,color=blue] (1,4)-- (2,3);
\draw [line width=1pt,color=blue] (1,2)-- (2,3);
\draw [line width=1pt,color=blue] (-3,4)-- (-4.04,2.84);
\draw [line width=1pt,color=blue] (-3,2)-- (-4.04,2.84);
\draw [fill=white] (-3,4) circle (0.2);
\draw [fill=white] (-1,4) circle (0.2);
\draw [fill=blue] (1,4) circle (0.2);
\draw [fill=blue] (-3,2) circle (0.2);
\draw [fill=blue] (-1,2) circle (0.2);
\draw [fill=blue] (1,2) circle (0.2);
\end{tikzpicture}
\end{center}
\caption{\label{fig3.23} Case $2.5$: $v_{i,j} \in S$, where $i \neq 1$, $2 < j < m-1$, and $j \equiv 2$ mod $3$}
\end{figure}
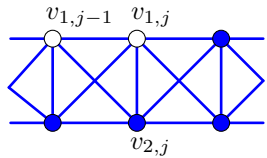
\indent Therefore, every vertex in $S$ has exactly two white neighbors. This proves that $S$ is stalled. Hence, $S$ is failed and $F(G) \geq |F|$ $=$ $ (nm - m) + \lceil \frac{m-4}{3} \rceil$. We claim that $S$ is maximal. \\
\indent Let $W$ be the set of vertices in the first copy of $P_m$. Hence, the induced subgraph $G[W]$ is a path with $m$ vertices. Since $v_{1,1}$ is an endpoint of the first copy of $P_m$, by Lemma \ref{lemma3.9}, the set $(W \cap S) \cup \{v_{1,1}\}$ is a standard zero forcing set for $G[W]$. Since all white vertices of $G$ are found in $W$, the set $S \cup \{v_{1,1}\}$ is a standard zero forcing set for $G$. Similarly, since $v_{1,m}$ is an endpoint of the first copy of $P_m$, we have that the set $S \cup \{v_{1,1}\}$ is a standard zero forcing set for $G$. Now, if $1 < j \leq 3\lceil \frac{m-4}{3}\rceil + 1$ for which $v_{1,j}$ is white, coloring $v_{1,j}$ blue gives rise to two consecutive blue vertices in the first copy of $P_m$. By Lemma \ref{lemma3.9}, the set $(W \cap S) \cup \{v_{1,j}\}$ is a standard zero forcing set for $G[W]$. Since all vertices of $G$ are in $W$, the set $S \cup \{v_{1,j}\}$ is a standard zero forcing set for $G$. It remains to show that $S \cup \{v_{1,j}\}$ is a standard zero forcing set for $G$ when $3\lceil \frac{m-4}{3}\rceil + 1 < j < m$. Let $S'$ be the set of vertices $v_{1,j}$ such that $3\lceil \frac{m-4}{3}\rceil + 1 < j < m$. We consider three cases. If $m-4$ $\equiv$ $0$ mod $3$, then $S' = \{v_{1,m-2}, v_{1,m-1}\}$. If $v_{1,m-2}$ is colored blue, then $v_{2,m-3}$ will force $v_{1,m-3}$. Now, $v_{1,m-4}$ and $v_{1,m-3}$ become two consecutive blue vertices in the first copy of $P_m$. By Lemma \ref{lemma3.9}, all the other white vertices will also be forced. If $v_{1,m-1}$ is colored blue, then $v_{2,m}$ will force $v_{1,m}$. Since $v_{1,m}$ is an endpoint, by Lemma \ref{lemma3.9}, all the other white vertices will also be forced. If $m-4$ $\equiv$ $1$ mod $3$, then $S' = \emptyset$. Finally, if $m-4$ $\equiv$ $2$ mod $3$, then $S' = \{v_{m-2}\}$. By a previous argument, coloring $v_{m-2}$ blue will allow forcing of all the other white vertices. \\
\indent Therefore, $S \cup \{v\}$ is a standard zero forcing set for $G$ for any vertex $v$ in $V(G) \backslash S$. Hence, $S$ is maximal.
\qed

\begin{theorem}\label{thm3.11} 
Let $G = C_m$ $\boxtimes$ $C_n$, where $m$ and $n$ are two integers greater than 2. Then, $F(G) \geq \lceil \frac{mn}{2} \rceil$.
\end{theorem} 
\proof
Let $G = C_n$ $\boxtimes$ $C_m$. The proof of Theorem \ref{thm3.5} gives us a failed zero forcing set $S$ for $C_n$ $\boxtimes$ $C_m$ with each vertex in $S$ having at least two neighbors in $C_n$ $\square$ $C_m$ outside $S$.  Hence, by Theorem \ref{thm3.8}, $S$ is also failed in $C_n$ $\boxtimes$ $C_m$ and $F(G) \geq |S| = \lceil \frac{mn}{2} \rceil$. 
\qed

\begin{theorem}\label{thm3.12} 
Let $n, m \in \mathbb{N}$ where $n, m \geq 2$. Then, $F(K_n$ $\boxtimes$ $K_m)$ $= nm-2$. 
\end{theorem}
\proof
We first show that $K_n$ $\boxtimes$ $K_m$ is just $K_{nm}$, the complete graph with $nm$ vertices. Suppose $G = K_n$ $\boxtimes$ $K_m$. Consider any pair of distinct vertices in $V(K_n) \times V(K_m)$, say $(u_1, v_1)$ and $(u_2, v_2)$. \\
\indent If $u_1 \neq u_2$ and $v_1 \neq v_2$, then, by definition of $K_n$ and $K_m$, $u_1$ is adjacent to $u_2$ and $v_1$ is adjacent to $v_2$. Hence, then, by definition of a strong product, $(u_1, v_1)$ and $(u_2, v_2)$ are adjacent. \\ 
\indent If  $u_1 = u_2$ and $v_1 \neq v_2$, then, by definition of $K_m$, $v_1$ is adjacent to $v_2$. Hence, then, by definition of a strong product, $(u_1, v_1)$ and $(u_2, v_2)$ are adjacent. \\
\indent Lastly, if  $u_1 \neq u_2$ and $v_1 = v_2$, then, by definition of $K_n$, $u_1$ is adjacent to $u_2$. Hence, then, by definition of a strong product, $(u_1, v_1)$ and $(u_2, v_2)$ are adjacent. \\
\indent Therefore, any two vertices in $G$ are adjacent. Hence, $G$ is a complete graph. Furthermore, $G = K_{nm}$ since there are $nm$ vertices in $K_n$ $\boxtimes$ $K_m$. Hence, by Theorem \ref{thm2.8}, $F(G) = F(K_{nm}) = nm - 2$. 
\qed

\setcounter{subsection}{2}
\subsection{Lexicographic products}
\setcounter{subsection}{12}

We next consider the failed zero forcing number of the lexicographic product of $G$ and $H$, denoted $G$ $\cdot$ $H$. Note that $G$ $\cdot$ $H$ and $H$ $\cdot$ $G$ are not necessarily the same graphs. 

\begin{theorem}\label{thm3.13} 
If $G$ and $H$ are graphs with $m$ and $n$ vertices, respectively, then \linebreak $F(G$ $\cdot$ $H)$ $\geq$ $s$ and there exists a set with cardinality $s$ that is a maximal failed zero forcing set for $F(G$ $\cdot$ $H)$, where
$$s = \begin{cases} 
nm-n+F(H), & \; if \; $H$ \; \text{has no isolated vertices} \\
nm-m+F(G), & \; if \; $H$ \; \text{has an isolated vertex} 
\end{cases}$$
\end{theorem}
\proof
Let $G$ and $H$ be graphs with $m$ and $n$ vertices, respectively. Consider the lexicographic product $G$ $\cdot$ $H$. \\
\indent Suppose $H$ has no isolated vertex. We construct a failed zero forcing set $F$ for \linebreak $G$ $\cdot$ $H$ as follows: Add to $F$ all vertices of $G$ $\cdot$ $H$ except the vertices in the $m$th copy of $H$. For the $m$th copy of $H$, include in $F$ only the vertices that form a maximum failed zero forcing set $F_H$ for $H$. Color the vertices in $F$ blue and the vertices in $V(G$ $\cdot$ $H)$ $\backslash$ $F$ white. Since $H$ has no isolated vertices, by Theorem \ref{thm2.7}, \linebreak $F(H) < n-1$. Hence, the derived coloring of $F_H$ in $H$ has at least two white vertices. Therefore, the product $G$ $\cdot$ $H$ initially has at least two white vertices in the $m$th copy of $H$ when $F$ is the initial set of blue vertices. \\
\indent Suppose $V(G) = \{a_1, ..., a_m\}$ and $V(H) = \{b_1, ..., b_n\}$. Then, $V(G$ $\cdot$ $H)$ $=$ $\{ (a_i, b_j)$ $:$ $1 \leq i \leq m$ and $1 \leq j \leq n\}$. Since $V(G$ $\cdot$ $H)$ $\backslash$ $F$ only contains vertices from the $m$th copy of $H$. Hence, a white vertex of $G$ $\cdot$ $H$ is of the form $(a_m, b_k)$, where $k \in \{1, ..., n\}$. Consider a blue vertex adjacent to $(a_m, b_k)$, say vertex $(a_r, b_s)$. It follows from the definition of a lexicographic product that $a_m$ is adjacent to $a_r$ in $G$, or $r = m$ and $b_s$ is adjacent to $b_k$ in $H$. \\
\indent If vertex $a_m$ is adjacent to vertex $a_r$ in $G$, then vertex $(a_r, b_s)$ is adjacent to vertex $(a_m, b_l)$ for $l \in \{1, ..., n\}$. Since the $m$th copy of $H$ has at least two white vertices, $(a_r, b_s)$ is adjacent to at least two white vertices of $G$ $\cdot$ $H$.  \\
\indent On the other hand, if $r = m$ and $b_s$ is adjacent to $b_k$, then $(a_r, b_s)$ is in the $m$th copy of $H$. Specifically, the vertex $(a_r, b_s)$ is in $F_H$. Since $F_H$ is maximum and failed, $(a_r, b_s)$ cannot force another vertex. Hence, $(a_m, b_k)$ cannot be forced. Since this holds for any $k \in \{1, ..., n\}$ such that $(a_m, b_k) \notin F$, $F$ is a failed zero forcing set for $G$ $\cdot$ $H$. Therefore, $F(G$ $\cdot$ $H)$ $\geq$ $nm - n + F(H)$. \\
\indent Now, consider $F \cup \{u\}$, where $u \notin F$. Then, $u$ is found in the $m$th copy of $H$ and $u \notin F_H$. Since all vertices except from the $m$th copy of $H$ are already blue and $F_H$ is already maximum, making $u$ blue allows forcing of the other white vertices found in the $m$th copy of $H$. Hence, $F \cup \{u\}$ is a standard zero forcing set for $G$ $\cdot$ $H$. Hence, $F$ is maximal. \\
\indent Next, suppose that $H$ has an isolated vertex. Then one of the connected components of $G$ $\cdot$ $H$, say $G_1$, has the same vertex set and edge set as $G$. Let $F$ be composed of all vertices in $V(G$ $\cdot$ $H)$ $\backslash$ $V(G_1)$ and all vertices of $G_1$ which forms a maximum failed zero forcing set $F_{G}$ for $G_1$. Color all the vertices in $F$ blue and all the vertices in $V(G$ $\cdot$ $H)$ $\backslash$ $F$ white. Consider a white vertex $v$ of $G_1$. Since $G_1$ is disconnected from the other components of $G$ $\cdot$ $H$, no vertex in \linebreak $V(G$ $\cdot$ $H)$ $\backslash V(G_1)$ can force $v$. Since $F_(G)$ is maximum and failed, no vertex in $V(G_1)$ can force $v$. Hence, vertex $v$ remains white. Therefore, $F$ is failed and $F(G$ $\cdot$ $H)$ $\geq$ $nm - m + F(G_1)$ $=$ $nm - m + F(G_1)$. Now, consider $F \cup \{u\}$, where $u \notin F$. Then, $u$ is found in the $n$th copy of $G$ and $u \notin F_G$. Since all vertices except from the $n$th copy of $G$ are already blue and $F_G$ is already maximum, making $u$ blue allows forcing of the other white vertices found in the $n$th copy of $G$. Hence, $F \cup \{u\}$ is a standard zero forcing set for $G$ $\cdot$ $H$. Hence, $F$ is maximal.
\qed

\indent The two corollaries follow immediately from Theorems \ref{thm2.8} and \ref{thm3.13}.

\begin{corollary}\label{cor3.14} 
Let $m, n \in \mathbb{N}$. Then, $F(P_m$ $\cdot$ $P_n)$ $\geq$  $nm - n + \lceil \frac{n-2}{2} \rceil$.
\end{corollary}

\begin{corollary}\label{cor3.15} 
Let $m$, $n \in \mathbb{N}$ $\backslash$ $\{1, 2\}$. Then, $F(C_m$ $\cdot$ $C_n)$ $\geq$ $nm - n + \lfloor \frac{n}{2} \rfloor$.
\end{corollary}

\begin{corollary}\label{cor3.16} 
Let $m$, $n \in \mathbb{N}$ $\backslash$ $\{1\}$. Then, $F(K_n$ $\cdot$ $K_m)$ $=$ $nm-2$. 
\end{corollary}
\proof
By Theorems \ref{thm2.8} and \ref{thm3.13}, $F(K_n$ $\cdot$ $K_m)$ $\geq$ $nm-n+(n-2)$ $=$ $nm-2$. Since $K_n$ and $K_m$ are both connected graphs, their lexicographic product is also connected. By Theorem \ref{thm2.7}, $F(K_n$ $\cdot$ $K_m)$ $<$ $|K_n$ $\cdot$ $K_m| - 1 $ $=$ $nm-1$. Hence, $F(K_n$ $\cdot$ $K_m)$ $=$ $nm-2$.
\qed

\setcounter{subsection}{3}
\subsection{Coronas}
\setcounter{subsection}{16}

\begin{theorem}\label{thm3.17} 
Suppose $G$ and $H$ are graphs with $m$ and $n$ vertices, respectively. If $W$ is the set of all isolated vertices of $H$ and $H'$ is the induced subgraph of $V(H) \backslash W$, then $F(G$ $\circ$ $H)$ $\geq$ $s$ and there exists a set with cardinality $s$ that is a maximal failed zero forcing set for $F(G$ $\circ$ $H)$, where
$$s = \begin{cases} 
2F(G), & \; if \; n=1 \\
nm+m-n+F(H), & \; if \; n \geq 2 \; \text{and} \; H \; \text{has no isolated vertices} \\
nm+m-2, & \; if \; n \geq 2 \; \text{and} \; H \; \text{is a graph with} \; n \; \text{isolated vertices} \\
nm+m-n+|W|+F(H'), & \; \text{otherwise}
\end{cases}$$
\end{theorem}
\proof
Let $G$ and $H$ be graphs with $m$ and $n$ vertices, respectively. Consider the corona $G$ $\circ$ $H$. \\
\indent Suppose $n=1$. Then $H = K_1$. Consider a maximum failed zero forcing set for $G$, say $F_G$. Let $F$ be composed of the vertices in $F_G$ and vertices in $F_H$, where $F_H$ consists of copies of vertices of $H$ corresponding to the vertices in $F_G$. Each vertex in $F_H$ is only adjacent to exactly one vertex in $F_G$. Hence, no vertex in $F_H$ can force another vertex. Since $F_G$ is failed and maximum, no vertex in $F_G$ can force another vertex. Hence, $F$ is failed. Consider the set $F \cup \{u\}$ for some $u \in V(G$ $\circ$ $H)$ $\backslash$ $F$. If $u$ is a vertex of $G$ which is not in $F_G$, then every vertex of $G$ will eventually be forced since $F_G$ is maximum and failed for $G$.  Since all vertices of $G$ are now forced and $H = K_1$, the vertex in each copy of $K_1$ will be forced by its corresponding vertex in $G$. On the other hand, suppose that $u$ is a copy of the vertex of $H$, then $u$ will force its corresponding vertex in $G$, say $v$. Since $v$ is not in $F_G$, by a previous argument, all vertices of $G$ $\circ$ $H$ will be forced. Therefore, $F \cup \{u\}$ is a standard zero forcing set for $G$ $\circ$ $H$. Thus, $F$ is maximal and $F(G$ $\circ$ $H)$ $\geq$ $2F(G)$. \\
\indent Suppose $n \geq 2$ and $H$ has no isolated vertex. We construct a failed zero forcing set $F$ for $G$ $\circ$ $H$ as follows: Add to $F$ all vertices of $G$ $\circ$ $H$ except the vertices in the $m$th copy of $H$. For the $m$th copy of $H$, include in $F$ only the vertices that form a maximum failed zero forcing set $F_H$ for $H$. Color the vertices in $F$ blue and the vertices in $V(G$ $\circ$ $H)$ $\backslash$ $F$ white. Since $H$ has no isolated vertices, by Theorem \ref{thm2.7}, $F(H) < n-1$. Hence, the derived coloring of $F_H$ in $H$ has at least two white vertices. Hence, if $F$ is the initial set of blue vertices, then $G$ $\circ$ $H$ initially has at least two white vertices in the $m$th copy of $H$. Consider a white vertex $v$ of $G$ $\circ$ $H$. It follows from the definition of a corona that if $u$ is the vertex of $G$ corresponding to $m$th copy of $H$, then no vertex in $G - u$ can force $v$. On the other hand, no vertex in the $m$th copy of $H$ can force $v$ since $F_H$ is maximum and failed. Hence, the white vertex $v$ cannot be forced. Since this holds for any vertex $v$ in $V(G$ $\circ$ $H)$ $\backslash$ $F$, $F$ is a failed zero forcing set for $G$ $\circ$ $H$. Therefore, $F(G$ $\circ$ $H)$ $\geq$ $nm + m - n + F(H)$. Now, consider $F \cup \{u\}$, where $u \notin F$. Then, $u$ is found in the $m$th copy of $H$ and $u \notin F_H$. Since all vertices except from the $m$th copy of $H$ are already blue and $F_H$ is already maximum, making $u$ blue allows forcing of the other white vertices found in the $m$th copy of $H$. Hence, $F \cup \{u\}$ is a standard zero forcing set for $G$ $\circ$ $H$. Hence, $F$ is maximal. \\
\indent Next, suppose $n \geq 2$ and $H$ is a graph with $n$ isolated vertices. Consider the vertices in the first copy of $H$. Since $n>1$, there are two vertices, say $u$ and $v$, in the said copy. Since every vertex in $H$ is isolated, it follows from the definition of a corona that $N(u) = N(v) = \{w\}$, where $w$ is the vertex of $G$ corresponding to the first copy of $H$. Hence, $\{u, v\}$ is a module of order 2. By Theorem \ref{thm3.5}, $F(G$ $\circ$ $H)$ $=$ $nm+m-2$. \\
\indent Finally, suppose that $n \geq 1$ and $H$ is a graph with at least one isolated vertex and at least one non-isolated vertex. Then $|W| > 1$, $|H'| > 1$, and $H'$ has no isolated vertex. Since $H'$ is a graph with more than one vertex, we can construct a maximal failed zero forcing set $F$ for $G$ $\circ$ $H'$ using the construction in the second case. Consider the set $F \cup W'$, where $W'$ are the vertices of $G$ $\circ$ $H'$ corresponding to the vertices in $W'$. Since the vertices in $W'$ are not adjacent to any vertex in $H'$, no vertex in $W'$ can force a white vertex in $H'$. By a previous argument in the second case, no vertex in $F$ can force a white vertex in $H'$. Hence, $F \cup W'$ is failed and $F(G$ $\circ$ $H)$ $\geq$ $m|H'| + m - |H'| + F(H') + |W'|$ $=$ $m|H'| + m - |H'| + F(H') + m(|H|-|H'|)$ $=$ $m|H| + m - (|H|-|W|) + F(H')$ $=$ $mn+m-m+|W|+F(H')$. Since $W'$ contains blue vertices disconnected from the white vertices of $G$ $\circ$ $H$, by the preceding paragraph, the set  $F \cup W' \cup \{u\}$ is a a standard zero forcing set for $G$ $\circ$ $H$ for any vertex $u$ in $V(G$ $\circ$ $H)$ $\backslash$ $(F \cup W')$. \\
\qed

\indent The next two corollaries then follow immediately from Theorems \ref{thm2.8} and \ref{thm3.17}.

\begin{corollary}\label{cor3.18} 
Let $m$, $n \in \mathbb{N}$. Then, $F(P_n$ $\circ$ $P_m)$ $\geq$ $nm + n - m + \lceil \frac{m-2}{2} \rceil$. 
\end{corollary}

\begin{corollary}\label{cor3.19} 
Let $m$, $n \in \mathbb{N}$ $\backslash$ $\{1, 2\}$. Then, $F(C_n$ $\circ$ $C_m)$ $\geq$ $nm + n - m + \lfloor \frac{m}{2} \rfloor$.
\end{corollary}

\begin{corollary}\label{cor3.20} 
Let $m$, $n \in \mathbb{N}$ $\backslash$ $\{1\}$. Then, $F(K_n$ $\circ$ $K_m)$ $=$ $nm + n -2$.
\end{corollary}
\proof
By Theorems \ref{thm2.8} and \ref{thm3.13}, $F(K_n$ $\circ$ $K_m)$ $\geq$ $nm+n-m+(m-2)$ $=$ $nm+n-2$. Since, $K_n$ and $K_m$ are both connected graphs, the corona $K_n$ $\circ$ $K_m$ is also connected. Hence, $F(K_n$ $\circ$ $K_m)$ $<$ $|K_n$ $\circ$ $K_m| - 1 $ $=$ $nm+n-1$. Hence, $F(K_n$ $\cdot$ $K_m)$ $=$ $nm+n-2$.
\qed

\setcounter{subsection}{4}
\subsection{Graphs showing that the bounds are sharp}
\setcounter{subsection}{20}

\begin{proposition}\label{prop3.21}
The bound in Theorem \ref{thm3.6} is sharp. 
\end{proposition}
\begin{enumerate}
\item Consider $G = P_2$ $\square$ $C_4$. By Theorem \ref{thm3.6}, $F(G) \geq \lceil \frac{4}{2} \rceil + 3 \lfloor \frac{4}{4} \rfloor$ $=$ $2+3$ $= 5$. Note that $G$ has no module of order 2 nor an isolated vertex. Since $G$ is connected, by Theorem \ref{thm2.7}, $F(G) < |G| - 2$ $=$ $2(4) - 2$ $=$ $6$. This implies that $F(G) \leq 5$. Therefore, $F(G) = 5$. 
\item Consider $G = P_2$ $\square$ $C_5$. By Theorem \ref{thm3.6}, $F(G) \geq \lceil \frac{5}{2} \rceil + 3 \lfloor \frac{5}{4} \rfloor$ $=$ $3+3$ $= 6$. A failed set with such cardinality is given by the set of blue vertices in the graph $G$ below. It can be verified that any set of seven vertices is a standard zero forcing set for the graph. Hence, $F(G) = 6$. 
\end{enumerate}

\begin{proposition}\label{prop3.22}
The bound in Corollary \ref{cor3.14} is sharp. 
\end{proposition}
\begin{enumerate}
\item By Corollary \ref{cor3.14}, $F(P_{10}$ $\cdot$ $P_{4})$ $\geq$ $40-4+1$ $=$ $37$. Observe that $P_{10}$ $\cdot$ $P_{4}$ is a connected graph and it has no module of order two, which implies that $F(P_{10}$ $\cdot$ $P_{4})$ $<$ $|P_{10}$ $\cdot$ $P_{4}| - 2$ $=$ $38$. Therefore, $F(P_{10}$ $\cdot$ $P_{4})$ $=$ $37$. 
\end{enumerate}

\begin{proposition}\label{prop3.23}
The bound in Corollary \ref{cor3.15} is sharp. 
\end{proposition}
\begin{enumerate}
\item By Corollary \ref{cor3.15}, $F(C_{3}$ $\cdot$ $C_4)$ $\geq$ $12-3+1$ $=$ $10$. Observe that $C_{3}$ $\cdot$ $C_4$ has no isolated vertices. Therefore, $F(C_{3}$ $\cdot$ $C_4)$ $<$ $|C_{3}$ $\cdot$ $C_4| - 1$ $=$ $12 - 1$ $=$ $11$. Therefore, $F(C_{3}$ $\cdot$ $C_4)$ $=$ $10$. 
\end{enumerate}

\begin{proposition}\label{prop3.24}
The bound in Corollary \ref{cor3.19} is sharp. 
\end{proposition}
\begin{enumerate}
\item By Corollary \ref{cor3.19}, $G$ $=$ $F(P_3$ $\circ$ $P_4)$ $\geq$ $12+3-4+1$ $=$ $12$. Note that $G$ is a connected and it has no module of order two, which implies that $F(G)$ $<$ $|G| - 2$ $=$ $15-2$ $=$ $13$. Therefore, $F(G)$ $=$ $12$. 
\end{enumerate}

\begin{proposition}\label{prop3.25}
The bound in Corollary \ref{cor3.20} is sharp. 
\end{proposition}
\begin{enumerate}
\item By Corollary \ref{cor3.20}, $G$ $=$ $F(C_4$ $\circ$ $C_3)$ $\geq$ $12+4-3+1$ $=$ $14$. Observe that $G$ is connected. Thus, it has no isolated vertices, which implies that $F(G)$ $<$ $|G| - 1$ $=$ $16-1$ $=$ $15$. Therefore, $F(G)$ $=$ $14$. 
\end{enumerate}

\section*{Acknowledgement}

The work of P. A. B. Pelayo was funded by the Accelerated Science and Technology Human
Resource Development Program - National Science Consortium (ASTHRDP-NSC) of the Department of
Science and Technology - Science Education Institute (DOST-SEI), Philippines.

\end{document}